\documentclass[10pt]{amsart}

\usepackage{amssymb, amscd, amsmath, amsthm,  graphicx, latexsym,enumerate}

\newtheorem{theorem}{Theorem}

\newtheorem{lemma}[theorem]{Lemma}

\newtheorem{corollary}[theorem]{Corollary}

 \theoremstyle{definition}
\newtheorem{definition}[theorem]{Definition}

\newtheorem{example}[theorem]{Example}

\begin{document}

\title[Relationships between braid length and the number of braid strands]
{Relationships between braid length and the \\number of braid strands}

\author{Cornelia A. Van Cott}
\address{Indiana University, Bloomington, Indiana, 47405} 
\email{cvancott@indiana.edu}
\maketitle


\section{Introduction}  

The first of many connections between knot theory and braid theory was discovered in 1923 when Alexander proved that every knot or link is isotopic to a closed braid~\cite{alex}.  The complicating twist is that this braid is not unique.  In fact, every knot is isotopic to infinitely many distinct closed braids.  

Many efforts have been made to study knots using their braid representatives.  In particular, braid theory is used to define two interesting invariants:  

\begin{definition}\label{index}The {\bf braid index} of a knot $K$, $br(K)$, is the minimum number of strands needed to express $K$ as a closed braid.  
\end{definition}

\begin{definition}\label{length}The {\bf braid length} of a knot $K$, $l(K)$, is the minimum number of crossings needed to express $K$ as a closed braid.  
\end{definition}

\begin{figure}\label{10-136}
 \includegraphics[width=60mm]{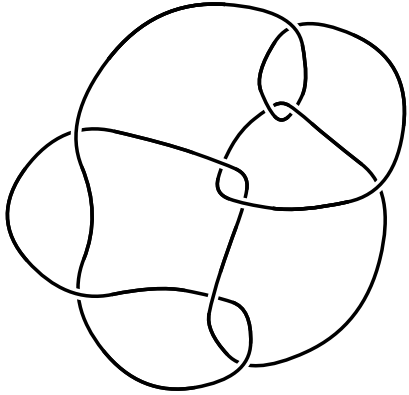}
 \caption{$10_{136}$:  Example of a knot with non-constant function $b_n'(K)$}
 \end{figure}

Consider the second knot invariant defined above:  braid length.  One might conjecture that a braid with the minimum number of crossings also has the minimum number of strands.  Thomas Gittings computed braid length for a large number of knots, and in the midst of these computations, he encountered counterexamples to this conjecture~\cite{gittings}.  The first counterexample is the knot $10_{136}$ (Figure~\ref{10-136}).  Gittings found that $br(10_{136})=4$ and $l(10_{136})=10$.  However, he observed that all 4--braids with closure $10_{136}$ have at least 11 crossings.  To find a 10--crossing braid with closure $10_{136}$, one must consider 5--braids.  The fact that adding strands might actually decrease the number of necessary crossings is unexpected, and in this paper we will investigate this balancing relationship between the number of braid strands and the length of braids associated to a knot $K$.  To this end, we give following definition:

\begin{definition}$b_n(K)=\text{minimum length of an $n$--braid representing }K$, where $n\geq br(K)$.  
\end{definition}


Fixing a knot $K$, we can consider $b_n(K)$ to be a function from a subset of the integers to the positive integers.
$$b_n(K):\bf {N}_{\geq br(K)}\longrightarrow\bf{N}$$
 We will study the behavior of this function.  Our interest will focus on knots for which the associated function decreases in value at some point, for these are precisely the knots for which increasing the number of strands actually decreases the number of necessary crossings.  The knot $10_{136}$ is the first known example of this.  From our previous discussion,
 
 \begin{eqnarray*}
b_4(10_{136})&=&11\\
b_5(10_{136})&=&10
\end{eqnarray*}

Theorem~\ref{difference} will show that $b_n(10_{136})=n+5$ for $n \geq5$.  So the function drops one time (from $b_4$ to $b_5$) and then increases monotonically thereafter.  As one might expect, knots such as $10_{136}$ for which the function $b_n(K)$ drops in value are anomalies among knots with low crossing number.  To ease the task of finding other examples similar to $10_{136}$, we consider the following fact.  If a knot $K$ is not one of the anomalies discussed above, then there exists some integer $a$ such that $b_n(K)=n + a$ for all $n \geq br(K)$ (this follows from Theorem~\ref{properties}).  So an arbitrary function $b_n(K)$ decreases in value at some point if and only if $b_n(K) -n$ is a non-constant function.  Studying $b_n(K)-n$ instead of $b_n(K)$ proves to be an easier way to think about our problem.  We define a new function as follows:


\begin{definition}$$b_n'(K):\bf {N}_{\geq br(K)}\longrightarrow\bf{N}$$
$$b_n'(K)=b_n(K)-n$$
\end{definition}

Whereas before we looked for functions $b_n(K)$ which dropped in value, we now are interested in finding functions $b_n'(K)$ which are not constant.  With these definitions in place, we state the following results to be proved in sections~\ref{firstSection} and~\ref{nextSection} about the behavior of $b_n'(K)$:

\begin{theorem}\label{properties}Let $K$ be a knot.  The function $b_n'(K)$ has the following properties:
\begin{enumerate}
 \item $b_n'(K)$ is a nonincreasing function.
 \item $b_n'(K)\equiv 1 \pmod{2}$.
 \item $b_n'(K) \geq 2g(K) -1$, where $g(K)$ = genus of the knot $K$.
  \end{enumerate}
\end{theorem}

\begin{theorem}\label{one}
Let $K$ be a nontrivial knot.  Then $b_n'(K)=1$ for some $n$ if and only if $K$ is the trefoil or the figure eight knot.  Hence, $b_n'(K) \geq 3$ for all $K$ not equal to the trefoil or figure eight.
\end{theorem}

\begin{theorem}\label{nine}
Let $K$ be a knot with nine or fewer crossings.  Then $b_n'(K)$ is a constant function.  
\end{theorem}

The final results, considered in Section~\ref{lastSection}, deal with two particular families of knots: homogeneous knots and braid positive knots~\cite{cromwell,stallings,s}.  Observe from the following definitions that the latter collection of knots lies inside the first.

\begin{definition}\label{homogeneous}
A knot $K$ is a {\bf homogeneous knot} if there exists a braid $\beta \in B_m$ for some $m$ with  its closure $\widehat{\beta}$ isotopic to $K$ such that given any braid group generator $\sigma_i$,  the exponent on $\sigma_i$ has the same sign in each appearance in the braid word $\beta$ (thus if $\sigma_i$ appears in the braid word $\beta$, then $\sigma_i^{-1}$ does not appear).  
\end{definition}

\begin{definition}\label{braidpositive}
A knot $K$ is a {\bf braid positive knot} if there exists a braid $\beta \in B_m$ for some $m$ with  $\widehat{\beta}$ isotopic to $K$ such that all generators of $B_m$ occurring in the braid word $\beta$ have positive exponents.
\end{definition}

 \begin{theorem}\label{homogeneousTheorem}
 Let $K$ be a homogeneous knot with an $m$--strand homogeneous braid representative.  Then $b_n'(K)=2g(K)-1$ for all $n \geq m$.
 \end{theorem}


\begin{theorem}\label{Stoimenow2}
Let $K$ be a braid positive knot with braid index $k$, such that $c(K) \leq 16$, where $c(K)$ denotes the crossing number of $K$.  If $b_k(K)\neq c(K)$, then $b_n'(K)$ is not a constant function. 
\end{theorem}

For background information about braids, see ~\cite{birman,birmanbrendle,murasugi}.  Except where otherwise noted, all knot enumeration and information is taken from~\cite{knotinfo}.

\subsection*{Acknowledgments} I thank Chuck Livingston for his assistance.  Thanks also to Mark Kidwell for his helpful comments and suggestions.  Jiho Kim wrote the computer program mentioned in Section 3. 
 
 \section{Properties of $b_n'(K)$}\label{firstSection}
To begin, we deduce several general properties of $b_n'(K)$ which aid in the computations of the function for specific knots.
\\
\\
\noindent{\bf Theorem~\ref{properties}. }{\it  Let $K$ be a knot.  The function $b_n'(K)$ has the following properties:
\begin{enumerate}
 \item $b_n'(K)$ is a nonincreasing function.
 \item $b_n'(K)\equiv 1 \pmod{2}$.
 \item $b_n'(K) \geq 2g(K) -1$, where $g(K)$ = genus of the knot $K$.
  \end{enumerate}}

\begin{proof}
Given a braid $\beta \in B_n$ whose closure is $K$, the stabilized braid $\beta \sigma_n \in B_{n+1}$ also has closure $K$.  Therefore, 

$$b_{n+1}(K) \leq b_n(K) +1.$$
Subtracting $n+1$ from both sides, we have $b_{n+1}'(K) \leq b_n'(K)$, proving $b_n'(K)$ is a nonincreasing function.

Any braid $\beta$ whose closure is $K$ determines a Seifert surface $S$ for the knot.  Suppose that the braid $\beta$ has $n$ strands and $c$ crossings.  Then the associated Seifert surface has $n$ Seifert circles and $c$ bands.  The genus of this surface is

$$g(S) = \frac{1+c-n}{2}.$$

Since the genus of a surface is an integer, it follows that $c-n\equiv 1 \pmod{2}$.  If $\beta \in B_n$ is a length minimizing $n$-braid for $K$, then $c=b_n(K)$.  It follows that $b_n'(K) = b_n(K)-n\equiv 1 \pmod{2}$.  Moreover, the Seifert surface $S$ determined by $\beta$ satisfies $g(K) \leq g(S) = \frac{1+b_n'(K)}{2}$.  Therefore, $b_n'(K) \geq 2g(K)-1$.

\end{proof}

\section{Computations of $b_n'(K)$}\label{nextSection}
Given the above basic properties of $b_n'(K)$, we can compute the function in specific cases.  For some knots, the function is straightforward to compute, as shown in the following lemma. 
 
\begin{lemma}\label{genus}
Let $K$ be a knot.  Let $\beta \in B_m$ be a braid of length $c$ whose closure is $K$.  Suppose $c-m=2g(K)-1$.  Then $b_n'(K)=2g(K)-1$ for all $n\geq m$.
\end{lemma}

\begin{proof}
By Theorem~\ref{properties}, $b_n'(K)$ is a nonincreasing function bounded below by $2g(K)-1$.  Since $b_m'(K)=2g(K)-1$, the conclusion follows.
\end{proof}

Though simple, this lemma applies to many knots, showing that their associated functions $b_n'(K)$ are constant.

\begin{example}\label{trefoil}
The trefoil $T$ is represented by the 2--braid $\alpha = \sigma_1^3$.  Applying Lemma~\ref{genus}, along with the fact that $g(T)=1$, it follows that $b_n'(T)=1$ for all $n\geq2$.  Thus the function $b_n'(T)$ is constant.  
\end{example}

\begin{example}\label{figure eight}
The figure eight knot $K$ can be represented by the 3--braid $\beta = \sigma_1\sigma_2^{-1}\sigma_1\sigma_2^{-1}$.  As with the trefoil, we can apply Lemma~\ref{genus} since $g(K)=1$ and conclude that $b_n'(K)=1$ for all $n\geq3$.  Since the figure eight knot has braid index $3$, the function $b_n'(K)$ is constant.

\end{example}

Notice that in both of the above examples, the function $b_n'(K)$ reaches the value 1.  Theorem~\ref{one} shows that these two knots, the trefoil and the figure eight, are the only knots for which this happens.
\\

\noindent{\bf Theorem~\ref{one}.} {\it  Let $K$ be a nontrivial knot.  Then $b_n'(K)=1$ for some $n$ if and only if $K$ is the trefoil or the figure eight knot.  Hence, $b_n'(K) \geq 3$ for all $K$ not equal to the trefoil or figure eight. }
\begin{proof}

From examples~\ref{trefoil} and ~\ref{figure eight}, the ``if" direction of the argument is clear. 

Conversely, let $K$ be a knot such that $b_n'(K)=1$ for some $n$.  Then $b_n(K)=n+1$ which implies that there is an $n$-braid $\beta$ with $n+1$ crossings and closure $\widehat{\beta}=K$.  

Count the number of occurrences of $\sigma_i^{\pm1}$ in the braid word $\beta$ for each $i\in\{1,2,\ldots,n-1\}$.  This gives us a string of $n-1$ integers which sum to $n+1$.  For example, the braid $\sigma_3\sigma_2\sigma_1\sigma_3^{-1}\sigma_4\sigma_2 \sigma_5^{-1}\in B_6$ corresponds to the string $(1,2,2,1,1)$.  

First, observe that none of the numbers in the string associated to $\beta$ are 0, for then $\widehat{\beta}$ would be a link.  By the pigeonhole principle, the string must be in one of the following forms:
$$(1,1,\ldots, 3,1\ldots,1)  \textrm{ or } (1,1,\ldots,2,\underbrace{1,\ldots,1}_{r},2,1,\ldots,1).$$  

If the string associated to $\beta$ is $(1,1,\ldots, 3,1,\ldots,1)$, then by conjugating and destabilizing the braid, we get a Markov equivalent $2$-braid $\beta'$ with $3$ crossings.  Since $\widehat{\beta}=\widehat{\beta'}=K$ is a non-trivial knot, it follows that $K$ is isotopic to the trefoil.

On the other hand, suppose the string is the latter of the two.  If $r>0$, then $K$ splits into the connected sum of two trivial knots and therefore is trivial, which is a contradiction.  So the string must be of the form $(1,1,\ldots,2,2,1,\ldots,1)$. By conjugating and destabilizing, we get a Markov equivalent braid with associated string $(2,2)$.  So $K$ can be expressed as a $3$-braid with $4$ crossings.  The only nontrivial knots which have diagrams with four crossings are the trefoil and the figure eight.  Therefore, $K$ must be either the trefoil or the figure eight.
 
Thus $b_n'(K)>1$ for all $n$, where $K$ is a non-trivial knot which is not the trefoil or figure eight.  By Theorem~\ref{properties} part 2, the function $b_n'(K)$ has odd integer values, so $b_n'(K)\geq 3$.

\end{proof}

\begin{example}
The knot $K=5_2$ has braid length $l(K)=6$ with braid representative $\beta=\sigma_1^3\sigma_2\sigma_1^{-1}\sigma_2 \in B_3$.  This means that $b_3(K)=6$ and $b_3'(K)=3$.  By Theorem~\ref{one}, $b_n'(K)=3$ for all $n\geq3$.  Since $K$ has braid index 3, we conclude that the function $b_n'(K)$ is constant.
\end{example}


In general, for any knot $K$ the function $b_n'(K)$ is eventually constant, since it is integer valued, nonincreasing, and bounded below.  The next result addresses the problem of finding an integer $n_0$ (depending on the knot $K$) such that the function $b_n'(K)$ is constant for all $n\geq n_0$.  The theorem finds such an integer $n_0$ in two different cases.  Notice that the two cases do not exhaust all possibilities;  however, if (using the notation of Theorem~\ref{difference}) $m = br(K)$, then the cases {\it are} exhaustive.


\begin{theorem}\label{difference}
Let $K$ be a prime knot.  Suppose $b_m(K)=c$ for some $m$.  Let $N = b_m'(K)=c-m$. 
\begin{enumerate}
\item If $N\geq m$, then $b_{n}'(K)=b_N'(K)$ for all $n \geq N$. 
\item If $N<br(K)$, then $b_{n}'(K)=b_m'(K)$ for all $n \geq m$. 
\end{enumerate}
\end{theorem}

\begin{proof}[Proof of (1)]
By assumption, there exists a braid $\beta \in B_m$ of length $c$ whose closure is isotopic to $K$.    

To begin, we show that $b_N(K)=2N-2t$ for some $t \geq 0$.  Notice that $b_m'(K)=c-m=N$ is odd by Theorem~\ref{properties} part 2.  Since $b_N'(K)=b_N(K)-N$ is also odd, it follows that $b_N(K)$ is an even integer.  Moreover, $b_N(K)$ is bounded above by $2N$, because we can take $\beta \in B_m$ and stabilize it $N-m$ times to get an $N$-braid whose closure is isotopic to $K$ with length $2N$.  Therefore $b_N(K)=2N-2t$ for some $t \geq0$.

Now we compute $b_{n}(K)$ where $n \geq N$.  Since $b_N(K)=2N-2t$, there exists a braid $\alpha \in B_{N}$ of length $2N-2t$ with closure isotopic to $K$.  As before, we can get an upper bound on $b_{n}(K)$ by taking $\alpha$ and stabilizing it $n-N$ times.  This produces a braid $\alpha' \in B_{n}$ of length $N-2t+n$ whose closure is isotopic to $K$.  Therefore, 
$$b_{n}(K)\leq N-2t+n.$$

Suppose that the above inequality is strict.  Then there exists a braid $\gamma \in B_{n}$ such that $\widehat{\gamma}$ is isotopic to $K$ and $l(\gamma)=N-2p+n$ where $p>t\geq 0$.  Now consider the string of $n-1$ positive integers associated to the braid $\gamma$, as in the proof of Theorem~\ref{one}.  The sum of the $n-1$ numbers in the string is $N-2p+n$, and none of the numbers are 0 (for otherwise $K$ would be a link).  A simple arithmetic calculation gives that the string has at least $(n-N)-2+2p$ entries equal to 1.  

The existence of a ``1" in the string implies that $\gamma$ must be the connected sum of two knots.  Since we assumed that $K$ is prime, it must be that one summand is the unknot.  So the braid is split by the single crossing into two sides:  one side of the braid corresponds to the unknot, and the other corresponds to $K$.  Since $\gamma$ has a minimal number of crossings, it follows that the part of the braid corresponding to the unknot also has a minimal number of crossings.  The unknot can be expressed as a $k$-braid via $\sigma_1\sigma_2\ldots\sigma_{k-1}$.  The integer string corresponding to this braid is $(1,1,\ldots,1)$, which implies that the braid has minimal length (a smaller length would force one of the integers to be 0).  Thus the integer string for $\gamma$ must be of the form $(1,1,\ldots,1,x_1,x_2,\ldots,x_r,1,\ldots,1)$, where none of the $x_i$'s  are 1.

The braid $\gamma$ can therefore be destabilized  $(n-N)-2+2p$ times.  In particular, we can destabilize the braid $n-N$ times, since $p>0$.  The resulting braid, $\gamma' \in B_N$ has length $2N-2p$.  Therefore, $b_N(K)\leq 2N-2p<2N-2t$.  This is a contradiction to $b_N(K)=2N-2t$, established earlier.  So we conclude 
\begin{eqnarray*}
b_n(K) & = & N-2t+n \\
& = & b_N(K)-N+n \\
& = & b_N'(K)+n
\end{eqnarray*}

Subtracting $n$ from both sides of this equation, it follows that $b_{n}'(K)=b_N'(K)$, as desired.
\\
\\
{\it Proof of (2).}  Suppose for contradiction that $b_n'(K) < b_m'(K)$ for some $n>m$.  This implies that $b_n'(K) = c-m-2p$ where $p>0$ and so $b_n(K) = c-m-2p+n$.  Hence, there is a braid $\beta \in B_n$ with $c-m-2p+n$ crossings whose closure is $K$.  Consider the associated string of $n-1$ positive integers for $\beta$.  As in Theorem~\ref{difference} part 1, the integers add up to $c-m-2p+n$, and since $K$ is a prime knot, none of the integers are zero and all of the 1's occur on either end of the string.  Again, a careful arithmetic calculation shows that the number of 1's in the string of integers is at least $n-c+m+2p-2$.  Thus, we can destabilize $\beta$ $n-c+m+2p-2$ times and obtain a braid $\beta'$ whose closure is $K$ and the number of strands in $\beta'$ is $c-m-2p+2 = N-2p+2 \leq N < br(K)$.  This contradicts the minimality of the braid index. 

\end{proof}

Theorem~\ref{difference} can be applied to the following special situation:  
\begin{corollary}\label{difference2}
Let $K$ be a prime knot of braid index $m$.  If $b_m(K)=l(K)$ and $b_m'(K)=m+1$, then $b_n'(K)=m+1$ for all $n\geq m$.  
\end{corollary}

\begin{proof}
We are given $b_m(K)=l(K)$.  Letting $N=b_m'(K)=m+1$, we can apply Theorem~\ref{difference} and conclude that $b_n'(K)=b_N'(K)$ for all $n \geq N=m+1$.  So to compute the function, it suffices to compute $b_N'(K)=b_{m+1}'(K)$.  We know that 
\begin{eqnarray*}
b_{m+1}'(K) &\leq& b_m'(K)\\
&=& b_m(K)-m\\
&=&l(K) -m.
\end{eqnarray*}
Suppose that the inequality is strict:  $b_{m+1}'(K) < l(K) -m$.  Then because of the parity of the function, $b_{m+1}'(K)<l(K)-(m+1)$, which implies $b_{m+1}(K) <l(K)$, a contradiction to the definition of braid length.  Therefore $b_{m+1}'(K)=b_m'(K)=m+1$.
\end{proof}

\begin{example}
Let $K=7_2$.  Though neither Theorem~\ref{one} nor Lemma~\ref{genus} apply to this knot, Corollary~\ref{difference2} does apply since $K$ has braid index 4 and $b_4(K)=9=l(K)$.  Therefore, $b_n'(K)=5$ for all $n\geq 4$. 
\end{example}

\begin{example}
In a manner similar to the previous examples, Theorem~\ref{one}, Lemma~\ref{genus}, Theorem~\ref{difference}, and Corollary~\ref{difference2} suffice to show that the function $b_n'(K)$ is constant for all except four knots with eight or fewer crossings.  The only outstanding knots are $7_3,7_5,8_{20},$ and $8_{21}$.  We can show that these four knots also have constant $b_n'(K)$ functions by using more computational methods as follows.  

First, we observe that these four knots have several identical knot invariants.  For $K\in \{7_3,7_5,8_{20},8_{21}\}$, we have $br(K)=3$, $l(K)=8$, $g(K)=2$, $b_3(K)=8$, and $b_4(K)=9$.  By Theorem~\ref{difference}, $b_5(K)$ determines $b_n(K)$ for all $n\geq5$.  So it remains to compute $b_5(K)$.  We know $b_5(K)$ is an even integer and has the following bounds: $8\leq b_5(K)\leq 10$.  Therefore, it suffices to show that $b_5(K)\neq 8$. 

Suppose that $b_5(K)=8$.  Then there exists a braid $\beta \in B_5$ of length 8 such that $\widehat{\beta}$ is isotopic to $K$.  Using an analysis similar to the proof of Theorem~\ref{one}, it follows that the string of 4 integers associated to $\beta$ is $(2,2,2,2)$.   

A computer check shows that none of the four knots can be formed by closing a braid whose string of integers is $(2,2,2,2)$.  This implies $b_5(K)\neq 8$, and therefore $b_5(K)=10$.  Now applying Theorem~\ref{difference}, it follows that $b_n'(K)=5$ for all $n\geq3$, where $K\in \{7_3,7_5,8_{20},8_{21}\}$.  So the function $b_n'(K)$ is constant for all knots $K$ with eight or fewer crossings.
\end{example}

\begin{example}
These results are sufficient for computing $b_n(K)$ and $b_n'(K)$ for nine crossing knots as well, though in this case almost half of the computations reduce to using the computer program mentioned above.  As before, the function $b_n'(K)$ is constant for {\sl all} nine crossing knots.  Hence for any knot $K$ with nine or fewer crossings, braid length cannot be decreased by increasing the number of braid strands.  
\end{example}

In view of this discussion, we have the following theorem:
\\
\\
\noindent{\bf Theorem~\ref{nine}.} {\it Let $K$ be a knot with nine or fewer crossings.  Then $b_n'(K)$ is a constant function.  }
\\
\\
The first known example of a non-constant function is $b_n'(10_{136})$, which was analyzed in the introduction.  Other knots for which the function is not constant include:  $11$n$_{8}, 11$n$_{121}, 11$n$_{131}, 12$n$_{20},12$n$_{24},12$n$_{65},12$n$_{119},12$n$_{358},12$n$_{362},$ and $12$n$_{403}$~\cite{knotinfo}.

\section{Homogeneous and braid positive knots}\label{lastSection}

We now turn our attention to two specific families of knots, homogeneous knots and braid positive knots (see definitions~\ref{homogeneous} and~\ref{braidpositive}), and to the proofs of theorems~\ref{homogeneousTheorem} and~\ref{Stoimenow2}.  Several special properties of these families enable us to find out more information about the associated functions $b_n'(K)$.  In particular, the following result of Stallings concerning homogeneous knots will be useful.  

\begin{theorem}[{\cite{stallings}}]\label{stallings}
If $\widehat{\beta}$ is the closure of a homogeneous braid, then $\widehat{\beta}$ is a fibred link.  
\end{theorem}


In the proof of the above theorem, Stallings also shows that the fibre surface $T_{\beta}$ for the knot $\widehat{\beta}$ is the Seifert surface corresponding to the closed braid $\widehat{\beta}$ with $m$ disks and $c$ bands, where $m$ (respectively, $c$) is the number of strands (crossings) of $\beta$.  Moreover, since $T_{\beta}$ is a fibre surface for $\widehat{\beta}$, it is a minimal genus Seifert surface:

\begin{equation}~\label{mingenus}
genus(T_{\beta}) = g( \widehat{\beta}).
\end{equation}

One further remark is in order before we begin the proof of Theorem~\ref{homogeneousTheorem}.  Given a homogeneous knot $K$ of braid index $k$, we cannot necessarily find a {\sl homogeneous} braid in $B_k$ whose closure is $K$.   Stoimenow found an infinite set of examples where no such braid exists ~\cite[Theorem 1]{s}.  In the definition, we are only guaranteed that for some integer $m \geq k$, there is a homogeneous braid in $B_m$ whose closure is $K$.  The following theorem shows that the function $b_n'(K)$ remains constant after this point $m$, and, in fact, we can state exactly what this constant function value is.
\\
 \\
\noindent{\bf Theorem~\ref{homogeneousTheorem}.} {\it    Let $K$ be a homogeneous knot with an $m$--strand homogeneous braid representative.  Then $b_n'(K)=2g(K)-1$ for all $n \geq m$.}

\begin{proof}
Let $l(\beta)=c$.  Let $S$ denote the Seifert surface (with $m$ disks and $c$ bands) corresponding to the closed homogeneous braid $\widehat{\beta}$.  Applying Theorem~\ref{stallings} and Equation~\eqref{mingenus}, we have
 $$g(K) = g(S)=\frac{1+c-m}{2}.$$
 
 Rearranging terms, we have $c-m=2g(K)-1$.  The following two inequalities follow from Theorem~\ref{properties} part 3 and the minimality of $b_m'(K)$:
 $$b_m'(K) \geq 2g(K)-1=c-m\geq b_m'(K).$$
 Therefore $b_m'(K)=2g(K)-1$.  Applying Lemma~\ref{genus}, we conclude $b_n'(K)=2g(K)-1$ for all $n\geq m$. 
\end{proof}

Our final result, Theorem~\ref{Stoimenow2}, gives some sufficient conditions for finding a knot $K$ such that $b_n'(K)$ is not constant.  First we state a result of Stoimenow:

\begin{theorem}[{\cite[Theorem 4]{s}}]\label{Stoimenow}
Any braid positive knot of crossing number $c(K)\leq 16$ has a braid positive minimal diagram.
\end{theorem}

The term {\sl braid positive minimal diagram} for a knot $K$ simply means a closed positive braid diagram in which the number of crossings equals $c(K)$.  Using this result, we prove Theorem~\ref{Stoimenow2}.
\\
\\
\noindent{\bf Theorem~\ref{Stoimenow2}.} {\it  Let $K$ be a braid positive knot with braid index $k$, such that $c(K) \leq 16$, where $c(K)$ denotes the crossing number of $K$.  If $b_k(K)\neq c(K)$, then $b_n'(K)$ is not a constant function.}

\begin{proof}
Since $b_k(K) \neq c(K)$, it follows that $b_k(K)>c(K)$.  By Theorem~\ref{Stoimenow}, there must be some integer $m> k$ and some positive braid $\beta \in B_m$ whose closure is $K$ such that $l(\beta)=c(K)$.  Hence $b_m(K) = c(K)$.  So we have
$$b_k(K) > b_m(K).$$
Subtracting $k$ from both sides and simplifying, we have:
\begin{eqnarray*}
b_k'(K)&=&b_k(K)-k \\
&>& b_m(K)-k  \\
&>& b_m(K)-m  \\
&=&b_m'(K).
\end{eqnarray*}
\end{proof}

\begin{figure}\label{p}
 \includegraphics[width=60mm]{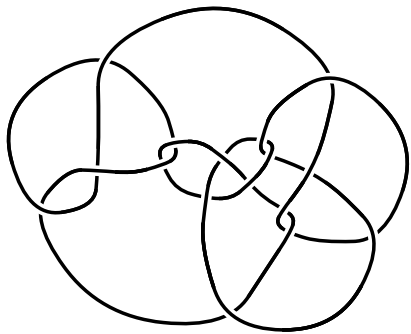}
 \caption{Braid positive knot $16_{472381}$}
 \end{figure}
 
\begin{example}
Consider the 16 crossing braid positive knot in Figure 2 ($K=16_{472381}$, see ~\cite{knotscape}).  The knot has braid index $4$, but no 4-braids in the equivalence class of $K$ are positive, and moreover all  have length at least $17$, which implies $b_4'(K)=13$.  Applying Theorem~\ref{Stoimenow2}, it follows that $b_n'(K)$ is not a constant function.  In fact, computations of ~\cite[Example 7]{s}  show that there is a positive 5-braid representing $K$ with length $16$.  So, $b_5'(K)=11$.  Now applying Theorem~\ref{homogeneousTheorem}, it follows that $b_n'(K)=11$ for all $n\geq 5$.  
\end{example}


\end{document}